\documentclass[letterpaper, 10 pt, conference]{ieeeconf}  

\IEEEoverridecommandlockouts                              

\usepackage[utf8]{inputenc}
\usepackage{bigints}
\usepackage{color,soul}
\usepackage{array}
\newcolumntype{P}[1]{>{\centering\arraybackslash}p{#1}}
\newcolumntype{M}[1]{>{\centering\arraybackslash}m{#1}}
\usepackage[thinlines]{easytable}
\usepackage{multirow}
\usepackage{xcolor,colortbl}

\definecolor{Gray}{gray}{0.85}
\definecolor{LightCyan}{rgb}{0.88,1,1}
\newcolumntype{a}{>{\columncolor{Gray}}c}
\usepackage{hyperref}
\usepackage{amsmath}
\usepackage{lipsum}
\usepackage{ams}
\usepackage{multicol, blindtext}
\newcommand{\norm}[1]{\left\lVert#1\right\rVert}

   \usepackage[pdftex]{graphicx}     
   \graphicspath{{../pdf/}{../jpeg/}{./image/}}    
   \DeclareGraphicsExtensions{.pdf,.jpeg,.png,.jpg}

\usepackage{amsmath}
\usepackage{bm}
\graphicspath{{Figures/}}
\usepackage[normalem]{ulem}
\usepackage[backend=biber,
style=ieee,
citestyle=numeric 
]{biblatex}
\graphicspath{ {./image/} }

\title{\LARGE \bf
Model and Load Predictive Control for Design and Energy Management of Shipboard Power Systems
}

\author{Mehrzad Mohammadi Bijaieh, Satish Vedula and Olugbenga Moses Anubi$^*$
\thanks{Authors are with the Department of Electrical and Computer Engineering, Center for Advanced Power systems, Florida State University, Tallahassee, FL 32310, USA
{\tt\small E-mail: \{mmohammadibijaieh, svedula, oanubi\}@fsu.edu}}%

}

\begin{document}

\maketitle
\thispagestyle{empty}
\pagestyle{empty}

\begin{abstract}
In current Medium Voltage DC (MVDC) Shipboard Power Systems (SPSs), multiple sources exist to supply power to a common dc bus. Conventionally, the power management of such systems is performed by controlling Power Generation Modules (PGMs) which include fuel operated generators and underlying converters. Moreover, energy management is performed by the emerging single or hybrid Energy Storage Systems (ESSs). This paper presents a model and load predictive control framework for power and energy management of SPSs. Here, MPC with load prediction is used for three main objectives: (1) to request power and energy from generators and energy storage elements according to their individual State of Power (SOP) and ramp-rate limitations, (2) to consider and integrate the generator cost and degradation, and (3) to reach a specific parking (final) State of Charge (SOC) for the ESSs at the end of the prediction horizon. The solution of the optimization problem is demonstrated using MATLAB and the functionality of the control framework is validated in real-time simulation environment.  
\end{abstract}

\section{INTRODUCTION}\label{Sec: Introduction}
Islanded MVDC Microgrids (MGs) deployed on SPSs provide power to various on-board equipment such as propulsion motors, hotel loads and highly non-linear loads such as Pulsed Power Loads (PPLs). There is a significant challenge to utilize existing ramp-rate limited generators to provide balanced power to high ramp-rate loads. The situation is exacerbated with integration of today's various power electronics equipment. This can result in an unbalanced system with degraded power quality, instability and subsequent load shedding and reconfiguration. 

A major challenge for operation of power electronics equipment is their Constant Power Load (CPL) behaviour where the current is inversely proportional to the voltage. In this case, a negative incremental impedance is created that can lead to instability  \cite{1976_MB,2009_Weaver}. However, fast utility load fluctuations cause more problems than existence of CPLs \cite{2015_Cupelli}. Since the ramp-rate support of generators are limited, ESSs with high ramp-rate support capabilities are offered as a solution \cite{2019_Bijaieh}.

Similar to DC MG systems, a hierarchical control architecture is currently utilized for power and energy management of SPSs \cite{2016_Jin}. In this case, a Power Management System (PMS) is designed to solve the underlying control allocation problem and the Energy Management System (EMS) is utilized for resource allocation. Conventionally, the PMS exists to ensure the stability and performance of the system while addressing the generation power over-actuation, and EMS exists to feed appropriate PMS set-points to achieve various high level objectives. This opens the ground for various optimization approaches for EMSs.   

Model Predictive Control (MPC) \cite{2009_Rawlings} is a mature technology \cite{2000_Mayne} that has been extensively used to control slow chemical processes. As the computation power improved throughout the years, it has become appropriate for control of faster systems such as power electronics. MPC utilizes the model of the system and its behaviour over a prediction horizon and aims to solve an optimization problem. It can systematically include the physical limitations such as State of Power (SOP) and ramp-rate limitations. Generally, MPC optimization includes an objective function with weighted expressions that often address a trade-off. The simplicity of addressing this trade-off comparing to the complexity of conventional algorithmic approaches is a major advantage of MPC control\cite{clemen2014model}. However, MPC problems tend to have computation burden that might not be viable in faster operations such as in real-time implementations\cite{anubi2015energy}. Hence, control system designers should always be aware of the computation cost of model predictive controllers.

Typically, there are different types of loads in a SPS. They can be categorized as linear, CPLs and high fluctuating utility loads. A characteristic of SPS is that some loads are rather cast than measured. The controllability of loads, the requirements of ESSs, and the importance of state awareness in such critical systems calls for integration of some form of predictive control. MPC-based EMSs for SPSs have been utilized in recent works such as \cite{2017_Vu}, where the ESSs ramp-rate limitation is readily programmed into the MPC problem. While a typical receding horizon MPC-based controller can achieve specific energy levels in ESSs, a planning approach for specific events and processes, as well as obtaining ESS requirements for sizing and technology specification, are viable paths to pursue.          

In this work, an MPC approach is used to demonstrate the generator-set (gen-set) and ESS control and resource allocation for a specific time horizon that can be as long as a mission operation or time duration specified for events such as undergoing a system-wide health-monitoring and degradation identification process through power excitation. In this work, the individual gen-set and battery SOP and ramp-rate, as well as a desired final SOC are considered as constraints for an objective function that aims to balance a trade-off between ESS processed energy and the gen-set operation cost. It will be shown that the MPC effectively allocates power to gen-set and the ESSs over a specific horizon and parks the corresponding SOC to a specific final value. The results of such control can be expanded and utilized by control designers to determine characteristics of individual system components such as generators and ESSs.

This work is organized as follows. Section \ref{Sec: System_Model} presents the SPS overall power flow and component simplified equivalent models. In Section \ref{Sec: Control_Development}, the MPC control as well as the component level control are shown. In Section \ref{Sec: Simulation} the proposed MPC with the corresponding objective function and constraints is first evaluated programmatically in Matlab and then its combination with the baseline component level control is implemented in Simulink Real-time.


\section{SYSTEM MODEL}\label{Sec: System_Model}
In this section, an example of a notional SPS \cite{esrdc1270} is demonstrated, the respective high-level equivalent model is presented and equivalent Low Bandwidth Model (LBM) of individual modules are specified.

\subsection{Notional Shipboard Power System}
\begin{figure}[t]
	\includegraphics[width=0.50\textwidth]{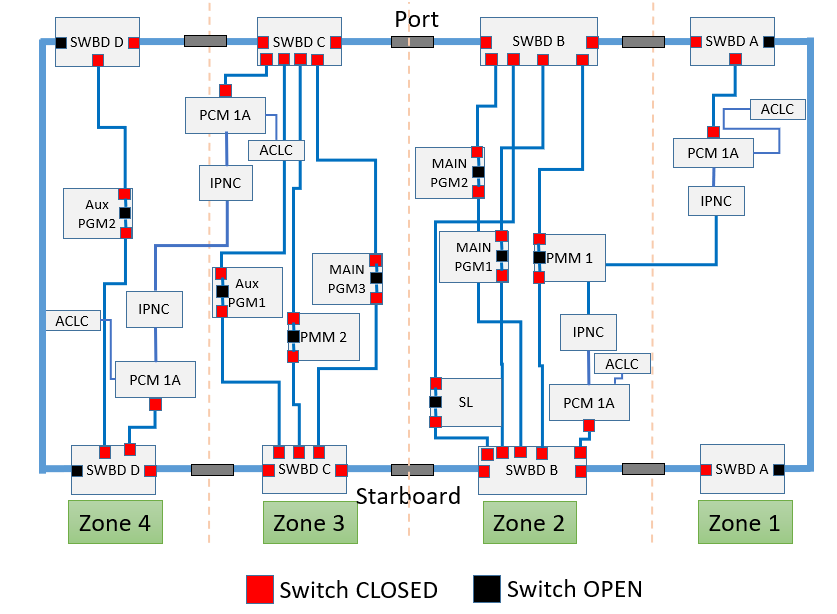}
	\caption{Notional 4-Zone SPS \cite{esrdc1270}}
    \label{4_Zone_SPS} \vspace{-3mm}
\end{figure}

SPS is considered as an islanded MVDC MG. Considering Fig. \ref{4_Zone_SPS}, the notional SPS is divided into four zones. Each zone may include one or a combination of power generating, storage and load modules. For example, zone-2 includes two main PGMs each including fuel operated generators as well as three-phase rectifiers, filters and respective Device Level Controllers (DLCs).   There are Power Conversion Modules (PCMs) which include multiple power converters, energy storage systems and AC and DC loads. High power and high ramp-rate loads are represented by Propulsion  Motor  Modules  (PMMs), Super Loads (SLs) and AC Load Centers (ACLC). Zones are connected to a unified $12kV$ DC power-line which should insure appropriate inter-zone power and energy transactions and zone and vessel-wide reconfiguration.


 \begin{figure}[h!] 
	\centering
	\includegraphics[width=0.35\textwidth]{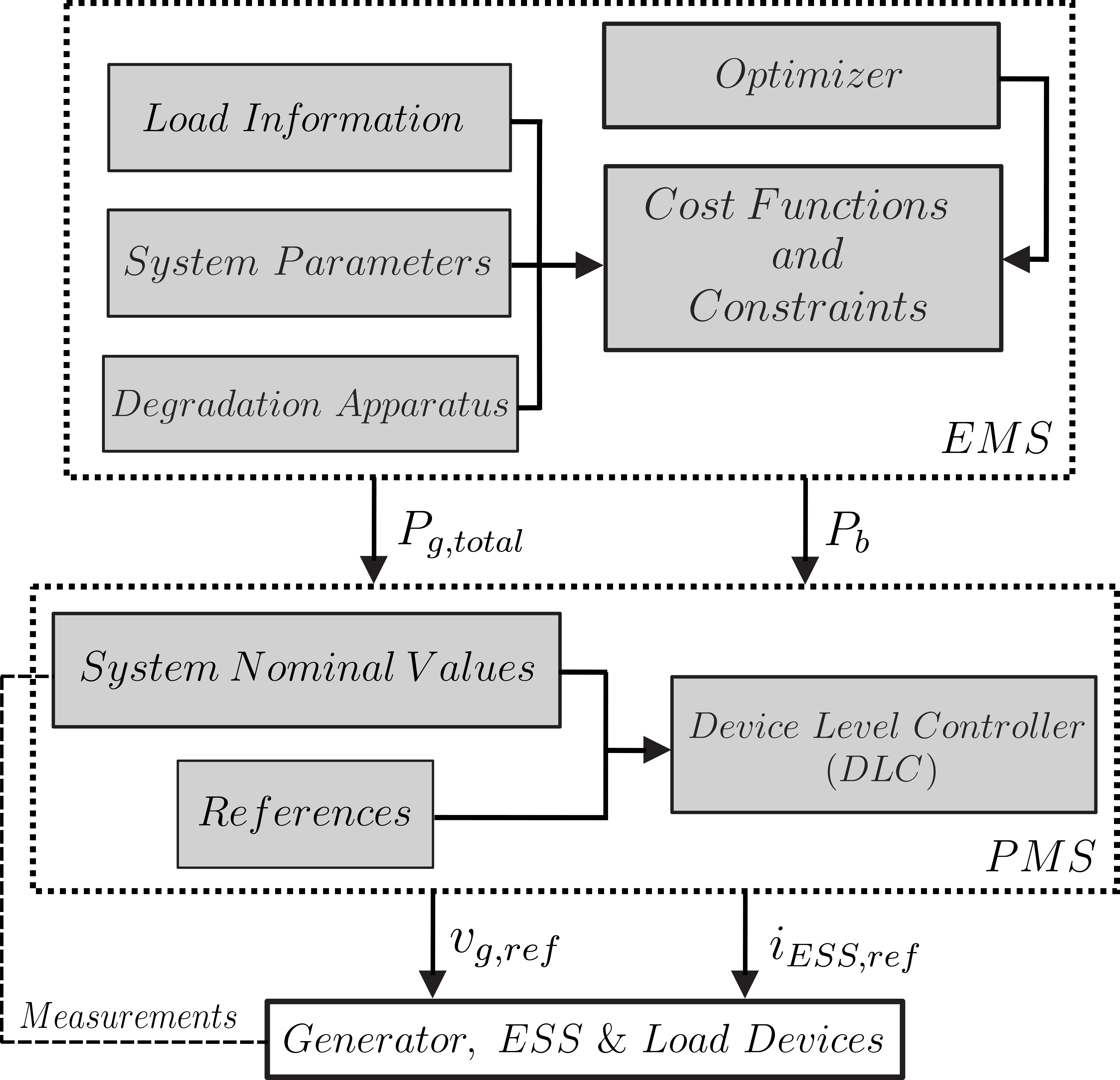}
	\caption{Control system diagram 
    }
	\label{Cntrl_Diag} \vspace{-3mm}
\end{figure}

\subsection{Power-flow Model}

The power-flow model must oblige to power generated equals power demanded criteria. Thus, The baseline SPS power-flow can be represented as:    
\begin{align} \label{SPS_P}
f(P_g,P_L,P_B)=0,
\end{align}
where, $P_g$, $P_B$ and $P_L$ represent the gen-set, batteries and load powers. In (\ref{SPS_P}), generators are considered to inject power in a unidirectional form while ESS can inject and absorb power bidirectionally to support the load demand. One major consideration for power distribution is ramp-rate limitation of generation and storage elements. Critical loads may have higher power ramp-rate demand than generators. A general approach is to request power from ESSs to compensate the high ramp-rate portion of the power demand which in turn enables the operation of the generators within their response capabilities. This is shown to improve the survivability, stability and quality of electrical levels of the SPS. Henceforth, the gen-sets and storage elements are referred to PGM and PCMs respectively.
\subsection{Component Model Representation}
\subsubsection{PGM}
The LBW model of the PGM includes a first-order filter fed into a controlled dependant voltage source, an $RL$ line, a shunt capacitor and a parallel damping resistor. The dynamic equations of PGM are as follows: 
\begin{subequations} \label{PGM_eqs}
\begin{align} 
    \frac{dv_g}{dt} &= k(v_{g,in}-v_g) \\ 
L_g\frac{di_g}{dt} &= v_g-i_gR_g-v_c  \\
C_g\frac{dv_c}{dt} &= i_g-\frac{v_c}{R_d} -i_{go},
\end{align}
\end{subequations}
where, $R_g$ and $L_g$  are the equivalent generator resistance and inductance respectively. $C_g$ is the generator capacitance, $v_g$ is generator voltage, $v_c$ is the voltage across the capacitor, $i_g$ is generator current, $R_d$ is the damping resistor, $k$ represents the cut-off frequency ($\frac{rad}{s}$) of the generator and $v_{g,in}$ denotes the lumped equivalent voltage control input of the PGM.

\subsubsection{PCM}
The PCM consists of an ideal ESS which may include single or hybrid storage systems such as Battery Energy Storage Systems (BESS), flywheels, or super-capacitor banks. The underlying equivalent dynamics for grid-following ESSs are given as:
\begin{align} \label{PCM_eqs}
\frac{{di}_{ESS}}{dt} &= \omega_{ESS}(i_{ESS,ref}-i_{ESS}),
\end{align}
where, $\omega_{ESS}$ is ESS response cut-off frequency. BESS is controlled through the reference current $i_{ESS,ref}$ and injects $i_{ESS}$ to the local bus. If the local bus voltage is shown as $v_b$, the reference ESS power reference can be shown as $P_{ESS,ref}=v_bi_{ESS,ref}$, while the injected power is $P_{ESS}=v_bi_{ESS}$ \cite{Bijaieh_2020}. 
The battery SOC calculations are given as:
\begin{subequations} \label{ESS_SOC}
\begin{align}
Q_{discharge} &= \frac{1}{3600} \int{i_{batt}(t)}dt \\
SOC &= \frac{Q_{0}-\frac{1}{3600}\int{i_{batt}(t)}dt}{Q_T},
\end{align}
\end{subequations}
where, $Q_{discharge}$ is the injected electric charge, and $Q_0$ and $Q_T$ are the the initial and total energy stored in the battery ESS in $AHr$. $i_{batt}$ is the injected battery current. The SOC versus the injected power can be obtained as
\begin{align} \label{SOCp}
SOC= \dfrac{Q_0v_b-\dfrac{1}{3600} \int P_{batt}(t)dt}{Q_Tv_b},
\end{align}
 where, $v_b$ represents the instantaneous measured voltage or the average voltage of the coupling bus.

\subsubsection{Load Module}
The simplified equivalent load is modelled as a controlled current source and a parallel $RC$ pair. The dynamics can then be specified as:
\begin{equation} \label{load_eq}
    C_L\frac{dv_{c,L}}{dt} = -\frac{P_i(t)}{v_{c,L}} - \frac{v_{c,L}}{{R_L}}+i_{in}
\end{equation}
where, $C_L$ is the load capacitance, $v_{c_L}$ is the load bus voltage, $P_{i}(t)$ is the load demand in terms of load power, $R_L$ is the resistive load, and $i_{in}$ is the overall load module current demand. 

For a PGM, PCM and cast load connections such as in Fig. \ref{Power_Flow}, the overall state-space representation for (\ref{PGM_eqs}), (\ref{PCM_eqs}) and (\ref{load_eq}) becomes
\begin{subequations} \label{SPS_all}
\begin{align} 
    \frac{dv_g}{dt} &= k(v_{g,in}-v_g) \\ 
L_g\frac{di_g}{dt} &= v_g-i_gR_g-v_{c,eq}  \\
\frac{{di}_{batt}}{dt} &= \omega_{batt}(i_{batt,ref}-i_{batt})\\
C_{eq}\frac{dv_{c,eq}}{dt} &= i_g+i_{batt}-\frac{v_{c,eq}}{R_d} -\frac{P_i(t)}{v_{c,eq}},
\end{align}
\end{subequations}
where, for a simplified series connection of the three modules, $C_{eq}$ is the equivalent shunt capacitance and is obtained from $C_{eq} = C_g + C_L$. The correspnding capacitor voltage is denoted by $v_{c,eq}$. The approach here is to control PGM to ensure voltage regulation at nominal value of $12kV$. Hence, the aim of the energy management through MPC is to control the ESS to portion the load with respect to its superior ramp-rate capabilities. Consequently, the rest of the load power is drawn from the bus which is regulated by PGM. The corresponding MPC and the baseline feed-forward and PI feedback control of PGM are shown in the next section.

\section{CONTROL DEVELOPMENT} \label{Sec: Control_Development}
The overall control hierarchy is shown in Fig. \ref{Cntrl_Diag}. MPC aims to optimize the SPS performance to meet a specific objective while considering the limitation or constraints that correspond to physical operability and performance. Hence, an objective cost function subject to system constraints is minimized over a specific stationary time horizon. 



\subsection{Model and Load Predictive Control}

The MPC objective function is 
	\begin{align} \label{MPC_obj}
\min_{P_{g,k},P_{b,k}}~~\sum_{k=1}^{h} \bigg(\norm{P_{g,k}+P_{b,k}-P_{L,k}^f}^2+\lambda C(P_{g,k})\bigg)
	\end{align}

The corresponding optimization constraints are 
\begin{subequations} \label{MPC_Const}	
	\begin{align} 
&\sum_{k=1}^{h}P_{b,k}=Q_b \\
&\underline{P_g} \leq P_{g,k} \leq \overline{P_g}\\
&| P_{g,k}-P_{g,k-1}| \leq r_g\\
&\underline{P_b} \leq P_{b,k} \leq \overline{P_b}\\
&| P_{b,k}-P_{b,k-1}| \leq r_b,
	\end{align}
\end{subequations} 
where, 
	\begin{align} \label{MPC_Const2}
Q_b= \dfrac{3600*Q_T*v_b^*}{T_s}(x_0-x_h). 
	\end{align}
where, $P_{g,k}$, $P_{b,k}$ and $P^f_{L,k}$ are the generator, battery and load forecast powers at $k^{th}$ horizon respectively. $C(P_{g,k})$ and $\lambda$ represent the cost of generator power and its weighting respectively. It is important to note that a large value for $\lambda$ increases the weight of generator cost in (\ref{MPC_obj}) which will eventually lead to a more active battery ESS. In the contrast, a zero value for $\lambda$ removes the generator cost considerations from the objective function. The MPC solves an optimization problem over the prediction horizon $h$ considering each time step denoted by $k$. The aim is to meet a specific load demand by using the gen-set and battery while considering generator operating cost. Here, the $C(P_{g,k})$ expression can be designed to incorporate generator or battery cost, efficiency or degradation. $C(P_{g,k})$ for the generator could also be related to the engine efficiency map or the risk associated in generating certain percentage of rated or maximum power. For the battery this could be effective degradation due to operating constantly at a certain C-Rate or the deviation from a given desired SOC. However, addressing the battery degradation and the subsequent ESS over-actuation is out of the scope of this work  and  is left for future iterations.	
Equation (\ref{MPC_Const2}) is obtained by discretizing (\ref{SOCp}). In (\ref{MPC_Const}), $\underline{P_g}$ and $\underline{P_b}$ represent minimum powers for generator and battery ESS respectively. Due to unidirectional operation of the generator, $\underline{P_g}$ cannot be less than zero while $\underline{P_b}$ can attain negative values since the battery operates bidirectionally. $\overline{P_g}$ and $\overline{P_b}$ denote the maximum injected powers respectively. 

(\ref{MPC_Const}a) is an equality constraint that represents the required injected power by the battery to reach a specific SOC through the prediction horizon $h$. (\ref{MPC_Const}b) and (\ref{MPC_Const}d) present the box constraints that correspond to SOP of the two systems. The equations (\ref{MPC_Const}c) and (\ref{MPC_Const}e) represent the absolute value of the difference between powers at current and previous step, and $r_g$ and $r_b$ denote the ramp-rate limit of the generator and the battery ESS. In (\ref{MPC_Const2}), $x_0$ and $x_h$ represent the initial and final SOC and $T_s$ is the corresponding time-step. 


The objective function (\ref{MPC_obj}) and the linear equality, rate and box constraints defined in (\ref{MPC_Const}) and (\ref{MPC_Const2}) can be directly reformed into the general quadratic form as 
\begin{align} \label{MPC_quad}
\min_{\mathbf{x}}~~\mathbf{x}^T \mathbf{H} \mathbf{x}+\mathbf{f}^T\mathbf{x},
\end{align}
such that 
\begin{subequations} \label{MPC_quad_const}
\begin{align} 
\mathbf{A_{eq}}.\mathbf{x}&=\mathbf{b_{eq}}\\
\mathbf{A}~.~\mathbf{x}&\preceq \mathbf{b}\\
\mathbf{x_{lb}}\preceq \mathbf{x} &\preceq \mathbf{x_{ub}},
\end{align}
\end{subequations}
and can be solved using quadratic programming \cite{2009_Rawlings}. It must be noted that MPC to QP formulation must comply with the above format. Many algorithms have been proposed mapping the MPC formulation into QP form \cite{6160293}.

\subsection{Baseline PGM Feedforward and Feedback Control}
The overall aim of PGM control is to regulate the main bus voltage to maintain the quality of electrical levels of the overall system and in extreme cases to ensure survivability of the SPS. The PGM feedforward control is obtained from the steady-state solution of the reference state-space system of (\ref{SPS_all}a), (\ref{SPS_all}b) and (\ref{SPS_all}d) such that:
\begin{subequations} \label{SPS_PGM_ref}
\begin{align} 
    \frac{dv_{g,ref}}{dt} &= k(v_{g,in}-v_{g, ref}) \\ 
L_g\frac{di_{g,ref}}{dt} &= v_{g,ref}-i_{g,ref}R_g-v_{c,eq}^*\\
C_{eq}\frac{dv_{c,eq}^*}{dt} &= i_{g,ref}-\frac{v_{c,eq}^*}{R_d} -\frac{P_i(t)}{v_{c,eq}^*}, 
\end{align}
\end{subequations}
where, the state index $``ref"$ denotes the corresponding reference state. Here, $v_{c,eq}^*$ is the nominal bus voltage value of $12kV$. Hence, the overall feedforward and feedback control for PGM is defined as
\begin{align} \label{SPS_PGM_Cntrl}
        v_{g,in} &= (\dfrac{R_g}{R_d}+1)v_{c,eq}^* +\dfrac{R_g}{v_{c,eq}^*}P_i(t)\nonumber\\
                 &~+ k_p * (v_{c,eq}^*-v_{c,eq})+ k_i\int_{0}^{\tau} (v_{c,eq}^*-v_{c,eq})d\tau,
\end{align}
where, $k_p$ and $k_i$ denote the proportional and integral gains respectively. $v_{c,eq}$ represents the measured load bus voltage.

\section{ILLUSTRATIVE EXAMPLES}\label{Sec: Simulation}
\begin{figure}[t!]
	\centering
 	\includegraphics[width=0.38\textwidth]{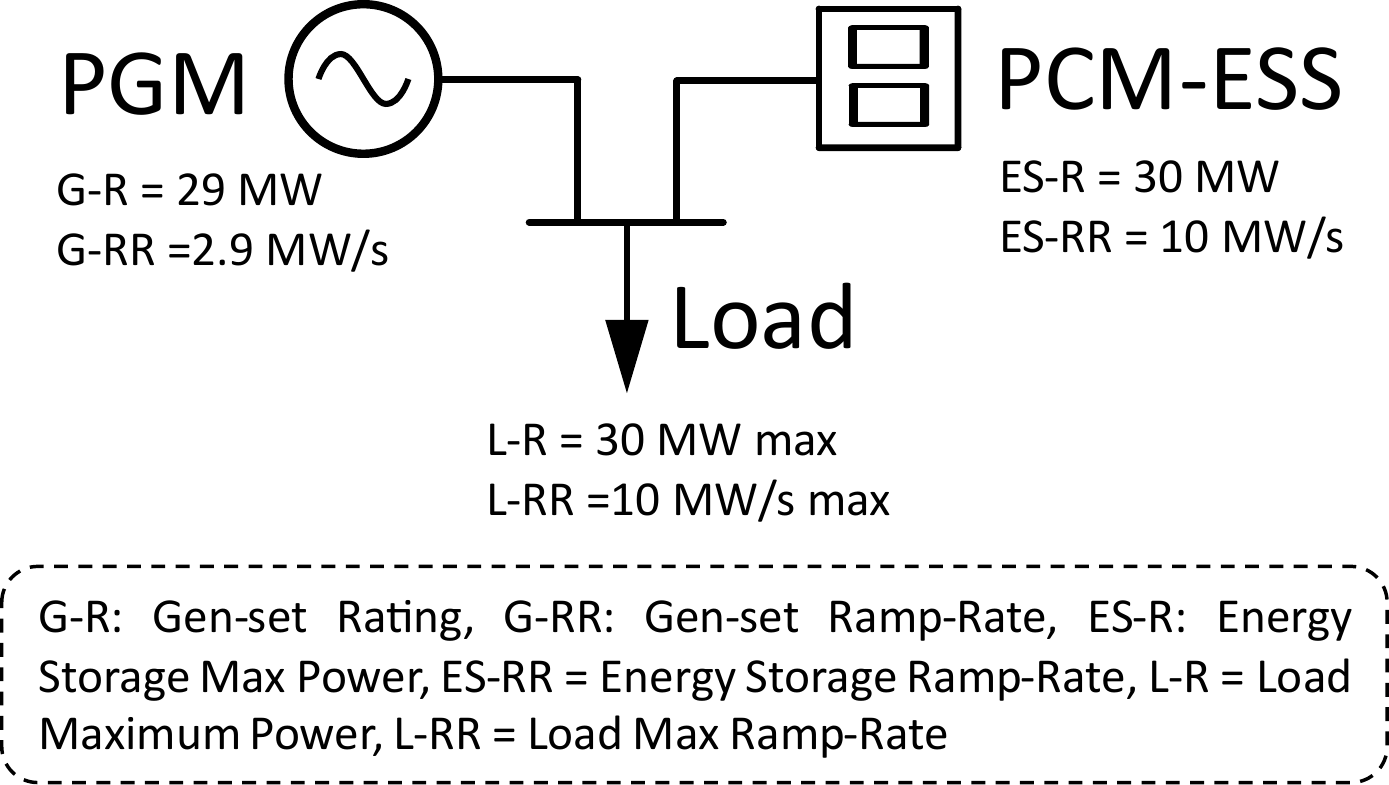}
 	\caption{Power system including generator set, a load and ESSs. 
     }
 	\label{Power_Flow} \vspace{-1mm}
 \end{figure}

The three-component system shown in Fig. \ref{Power_Flow} is used to demonstrate the solution of optimization problem defined in (\ref{MPC_obj}), (\ref{MPC_Const}) and (\ref{MPC_Const2}) over a specific time horizon. The reformed MPC problem in (\ref{MPC_quad}) with constraints in (\ref{MPC_quad_const}) is first evaluated programmatically in MATLAB, then, the overall control and resource allocation behavior is implemented in a Software-in-the-Loop (SIL) environment on a SN5726 Speedgoat performance real-time target. The corresponding real-time data is obtained and then plotted using MATLAB. The rated powers as well as the ramp-rate specifications for the generator, battery and the load are shown in Fig. \ref{Power_Flow}. The PGM, PCM and load system models are defined in (\ref{SPS_all}) with the baseline DLC shown in (\ref{SPS_PGM_Cntrl}). The corresponding system and control parameters are shown in Table \ref{sys_Param}. 

\begin{figure}[t!] 
 	\centering
 	\includegraphics[width=0.46\textwidth]{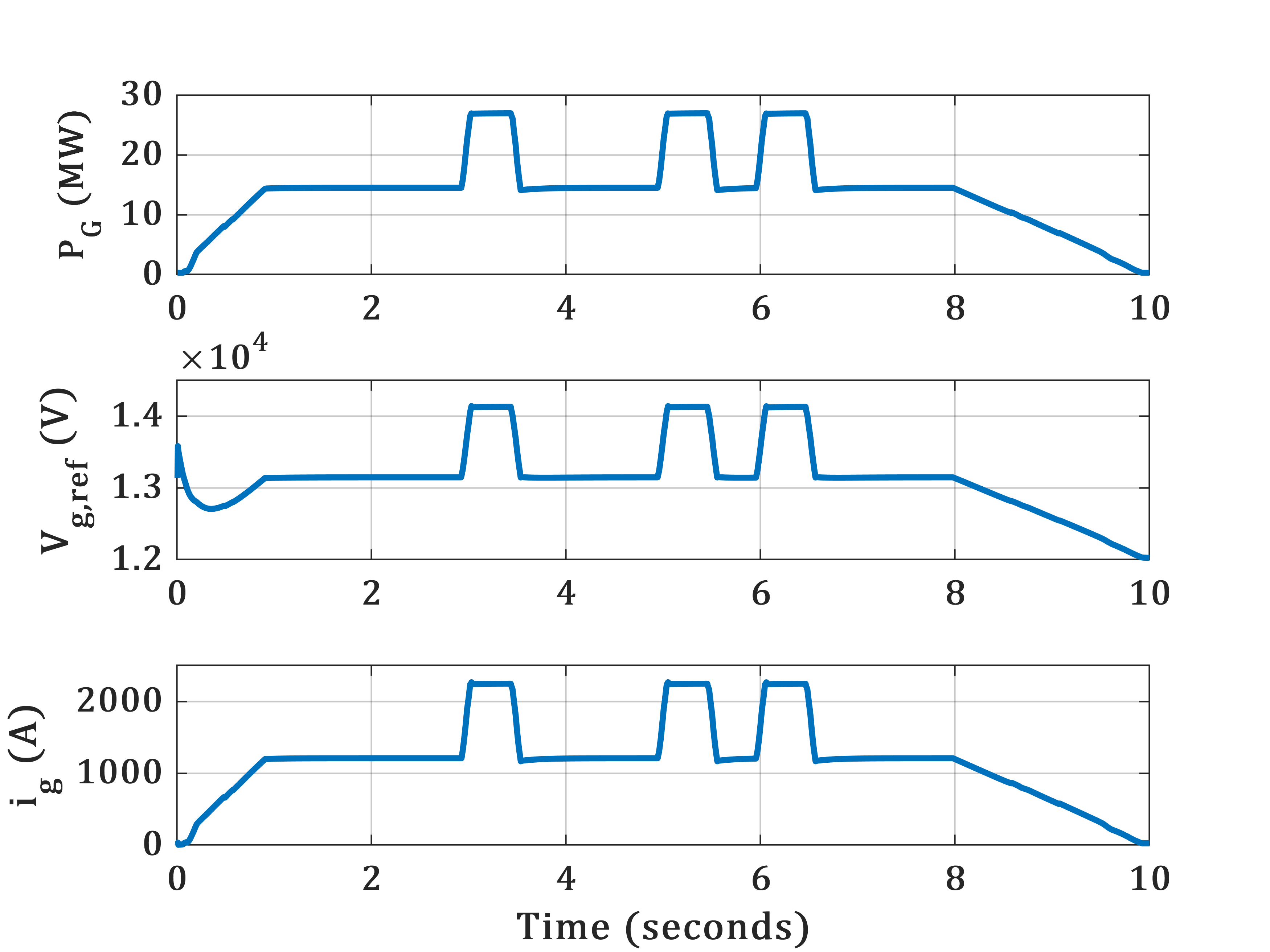}\vspace{-3mm}
 	\caption{(a) Power injected by the gen-set (b) Generator voltage reference command (c) Generator real-time current for initial SOC of 0.8 and final SOC of 0.77
     }
 	\label{Generator08-077} \vspace{-4mm}
 \end{figure}

\begin{figure}[t!]
 	\centering
 	\includegraphics[width=0.48\textwidth]{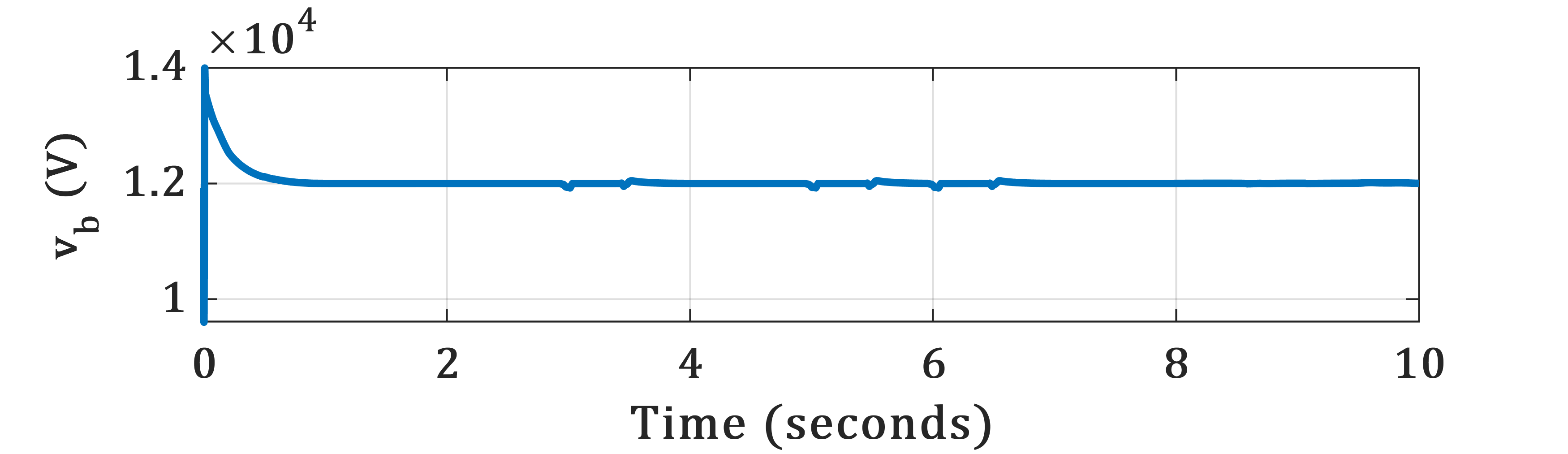} \vspace{-3mm}
 	\caption{ Load bus voltage in real-time controlled around 12kV 
     }
 	\label{Vb08-077} \vspace{-3mm}
 \end{figure}

\begin{figure}[t!]
 	\centering
 	\includegraphics[width=0.46\textwidth]{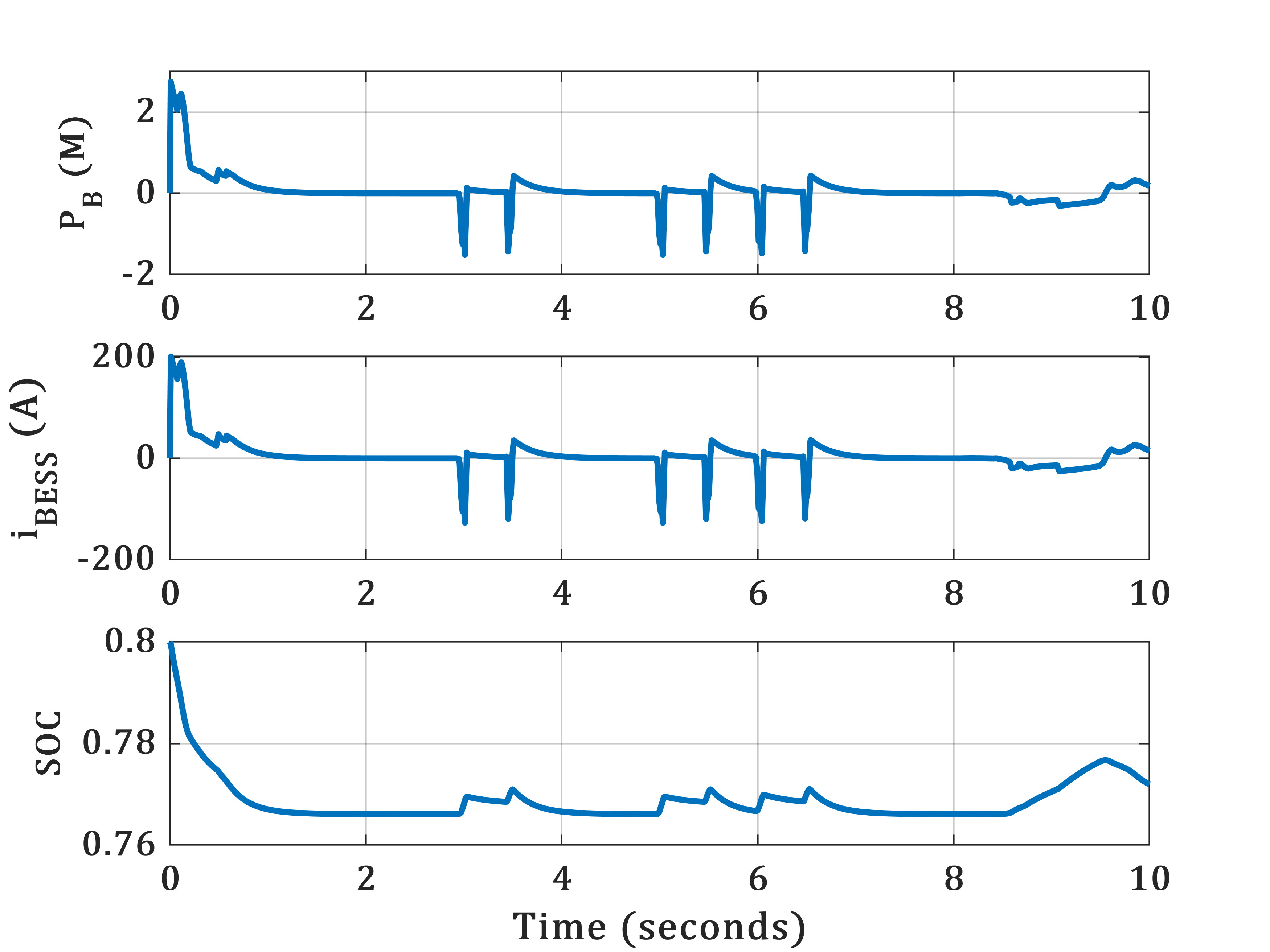} \vspace{-3mm}
 	\caption{(a) Power injected by the ESS (b) ESS real-time current (c) SOC variation for initial SOC of 0.8 and final SOC of 0.77
     }
 	\label{BESS08-077} \vspace{-3mm}
 \end{figure}

\begin{figure}[t!]
 	\centering
 	\includegraphics[width=0.46\textwidth]{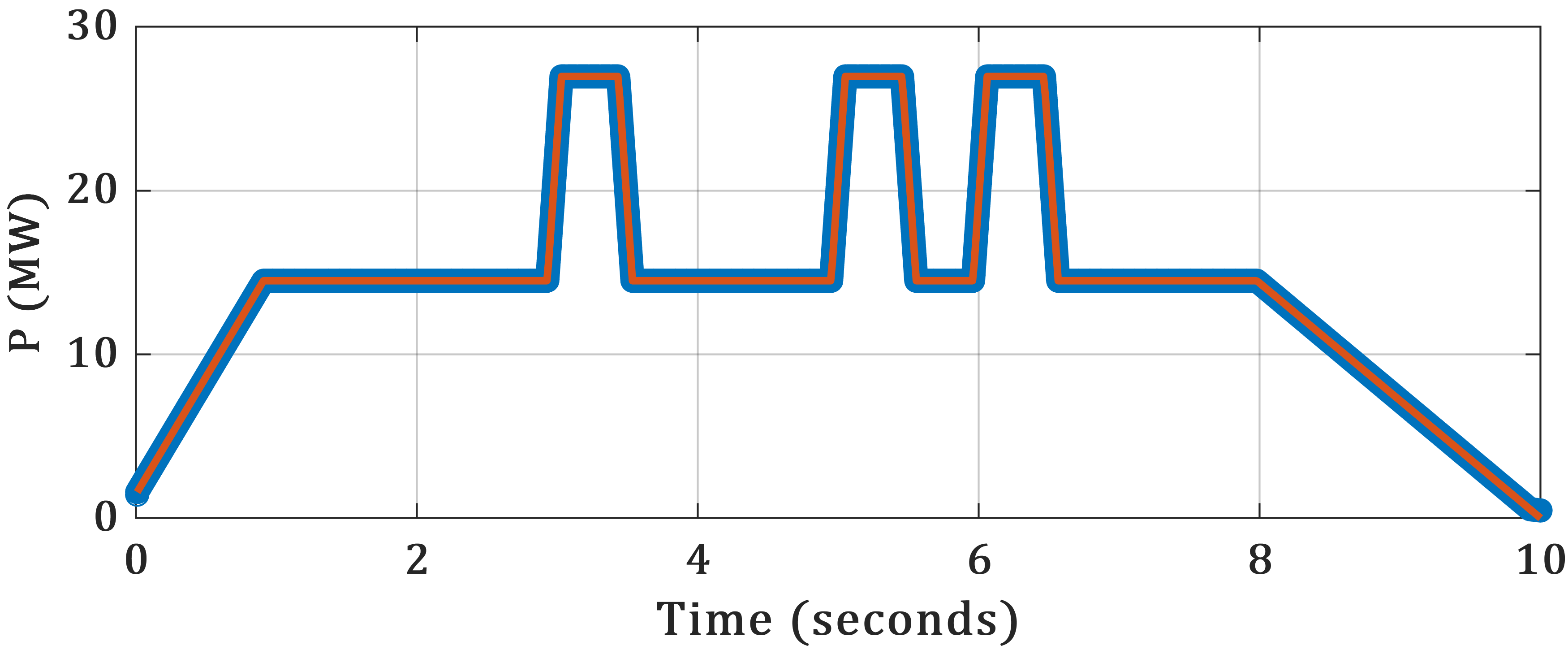} \vspace{-3mm}
 	\caption{ Power injected by the combined gen-sets and ESS matches the power forecast 
     }
 	\label{PD0808} \vspace{-3mm}
 \end{figure}

The over all results are shown for two cases: (1) To maintain SOC of the battery system and keep the parking SOC of the battery system at 0.8, and (2) to start with a 0.8 SOC and reach a specific parking SOC of 0.77. The results for the first and second cases are demonstrated in Figs. \ref{Generator08-08} to \ref{BESS0808} and Figs. \ref{Generator08-077} to \ref{BESS08-077} respectively. In both cases same load profile is used for critical comparison. Here, load profile is composed of superimposed signals that mimic hotel and high ramp-rate controllable loads such as PPLs. For both cases, one objective is to utilize the ESS when the ramp-rate of the generator is not sufficient to match the load. Also, observe the behavior of ESS for changes in parking SOC.

    \begin{table}[t]
\renewcommand{\arraystretch}{1.3}
\centering
\caption{System \& Control Parameters} \vspace{-3mm}
\label{sys_Param}
\begin{tabular}{lcccccccc}
\hline
\hline
\rowcolor{Gray}
\multicolumn{9}{c}{\textbf{PGM, PCM and Load Parameters}}                                                 \\ \hline 
               & \multicolumn{2}{c}{$R_g(\Omega)$}  & \multicolumn{2}{c}{$R_d(k\Omega)$}  & \multicolumn{2}{c}{$L_{g}(H)$} & \multicolumn{2}{c}{$k$}  \\ \hline
 PGM    & \multicolumn{2}{c}{0.9}  & \multicolumn{2}{c}{10}  & \multicolumn{2}{c}{0.1} & \multicolumn{2}{c}{10}  \\ 
 \hline
 & \multicolumn{2}{c}{$\omega_{batt}(\frac{rad}{s})$}  & \multicolumn{2}{c}{$Q_{T}(GJ)$}  & \multicolumn{2}{c}{$C_{g} (m F)$} & \multicolumn{1}{c}{$~$}  \\ \hline
 PCM \& Load    & \multicolumn{2}{c}{1000}  & \multicolumn{2}{c}{20}  & \multicolumn{2}{c}{2.1} & \multicolumn{1}{c}{~}  \\ \hline
 
 \rowcolor{Gray}
 \multicolumn{9}{c}{\textbf{Control Parameters}}                                                      \\ \hline
               & \multicolumn{2}{c}{$k_p$}  & \multicolumn{2}{c}{$k_i$}  & \multicolumn{2}{c}{$v_{c,eq}^*(kV)$} & \multicolumn{2}{c}{$~$}  \\ \hline
 PCM \& Load    & \multicolumn{2}{c}{1}  & \multicolumn{2}{c}{10}  & \multicolumn{2}{c}{12} & \multicolumn{2}{c}{~}  \\ \hline \hline

\end{tabular}
\end{table}  


\subsection{Parking SOC as Initial SOC}

\begin{figure}[t!]
 	\centering
 	\includegraphics[width=0.46\textwidth]{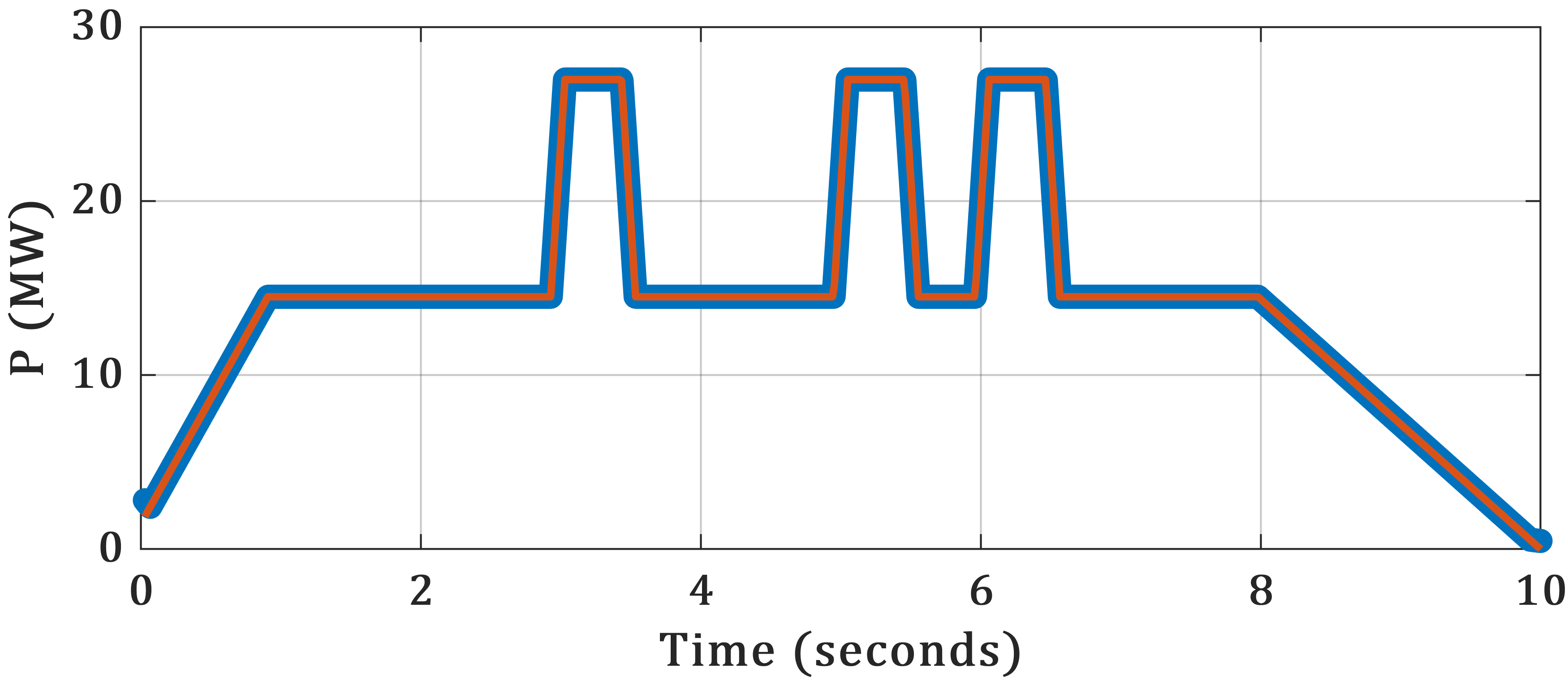} \vspace{-3mm}
 	\caption{ Power injected by the combined gen-sets and ESS matches the power forecast 
     }
 	\label{PD08-077} \vspace{-3mm}
 \end{figure}

In this case, ESS is controlled to inject a zero average power at the end the prediction horizon. For this study the following simulation parameters are considered: horizon $h=5$, time step for simulation is $100\mu s$. The time scale separation between energy management layer the system layer is considered to be $1ms$, this accounts for device level component time constants. As shown in Fig. \ref{PD08-077}, load includes fast fluctuations as well as hotel loads. Both the generator and the ESS contribute to match the load power. Here, the generator injects nearly all of the total power except when it cannot keep up with the high ramp-rate demand. On the other hand, ESS injects most of its power when there is high ramp-rate demand. 

\begin{figure}[t!]
 	\centering
 	\includegraphics[width=0.46\textwidth]{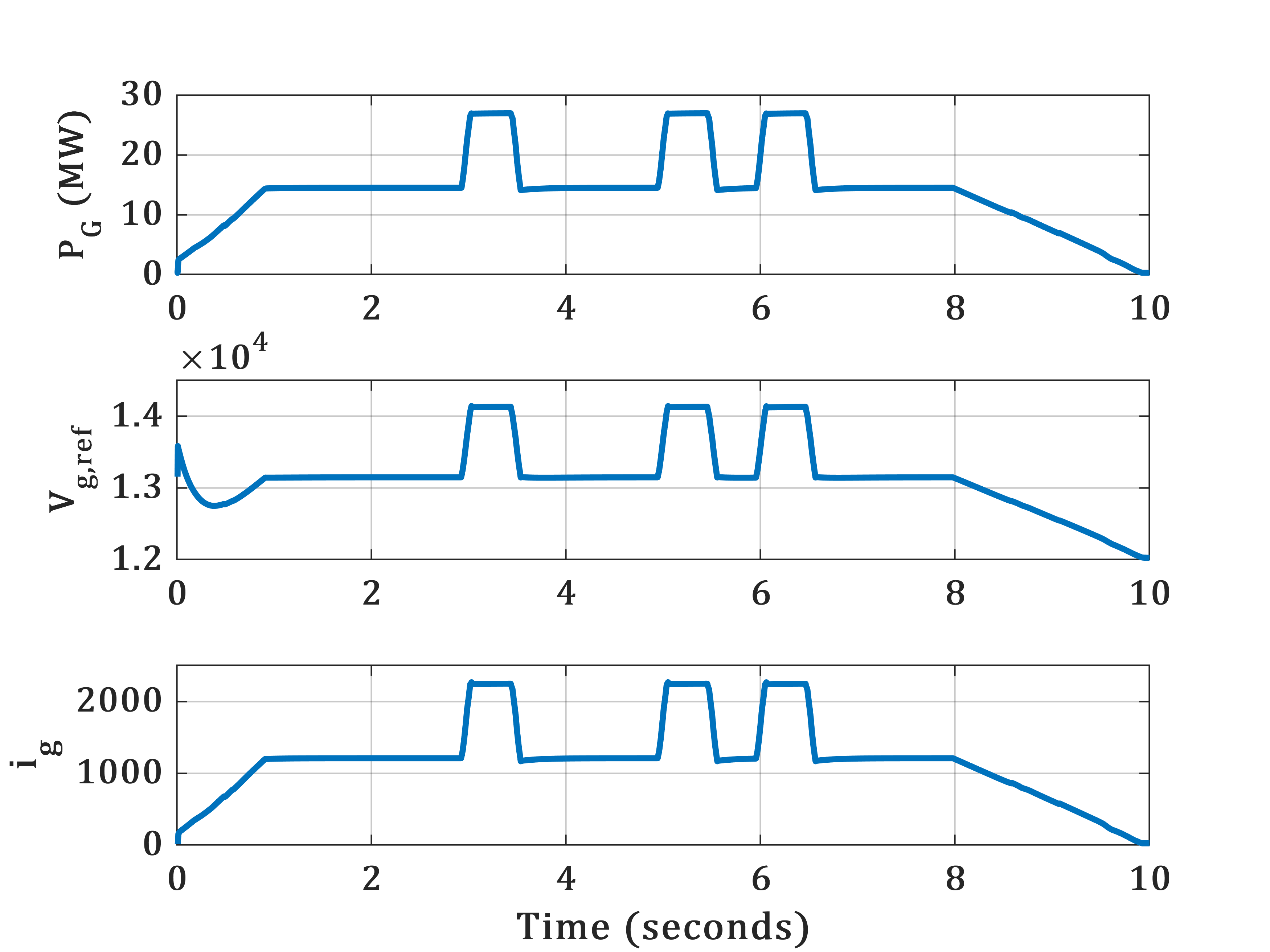} \vspace{-3mm}
 	\caption{(a) Power injected by the gen-set (b) Generator voltage reference command (c) Generator real-time current for initial SOC of 0.8 and final SOC of 0.8 
     }
 	\label{Generator08-08} \vspace{-4mm}
\end{figure}
\begin{figure}[t!]
 	\centering
 	\includegraphics[width=0.46\textwidth]{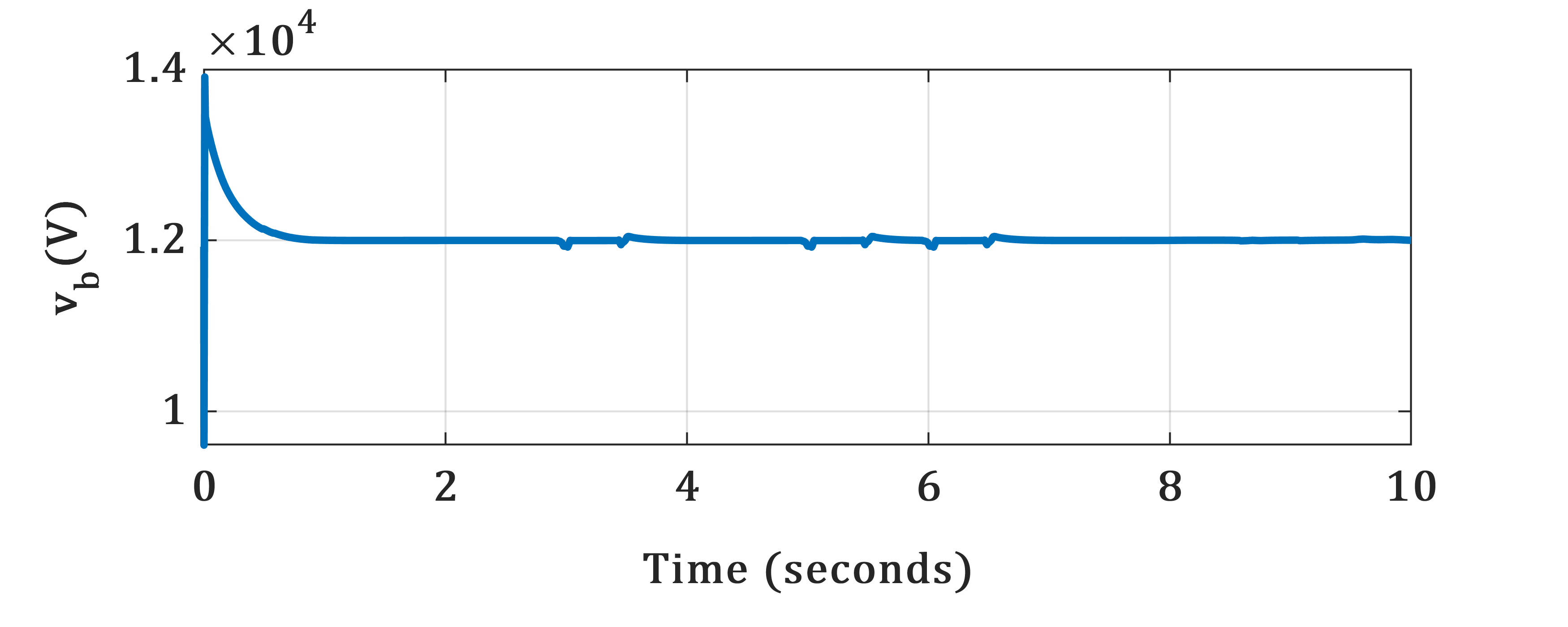} \vspace{-3mm}
 	\caption{ Load bus voltage in real-time controlled around 12kV 
     }
 	\label{Vb0808} \vspace{-3mm}
 \end{figure}
\begin{figure}[t!]
 	\centering
 	\includegraphics[width=0.46\textwidth]{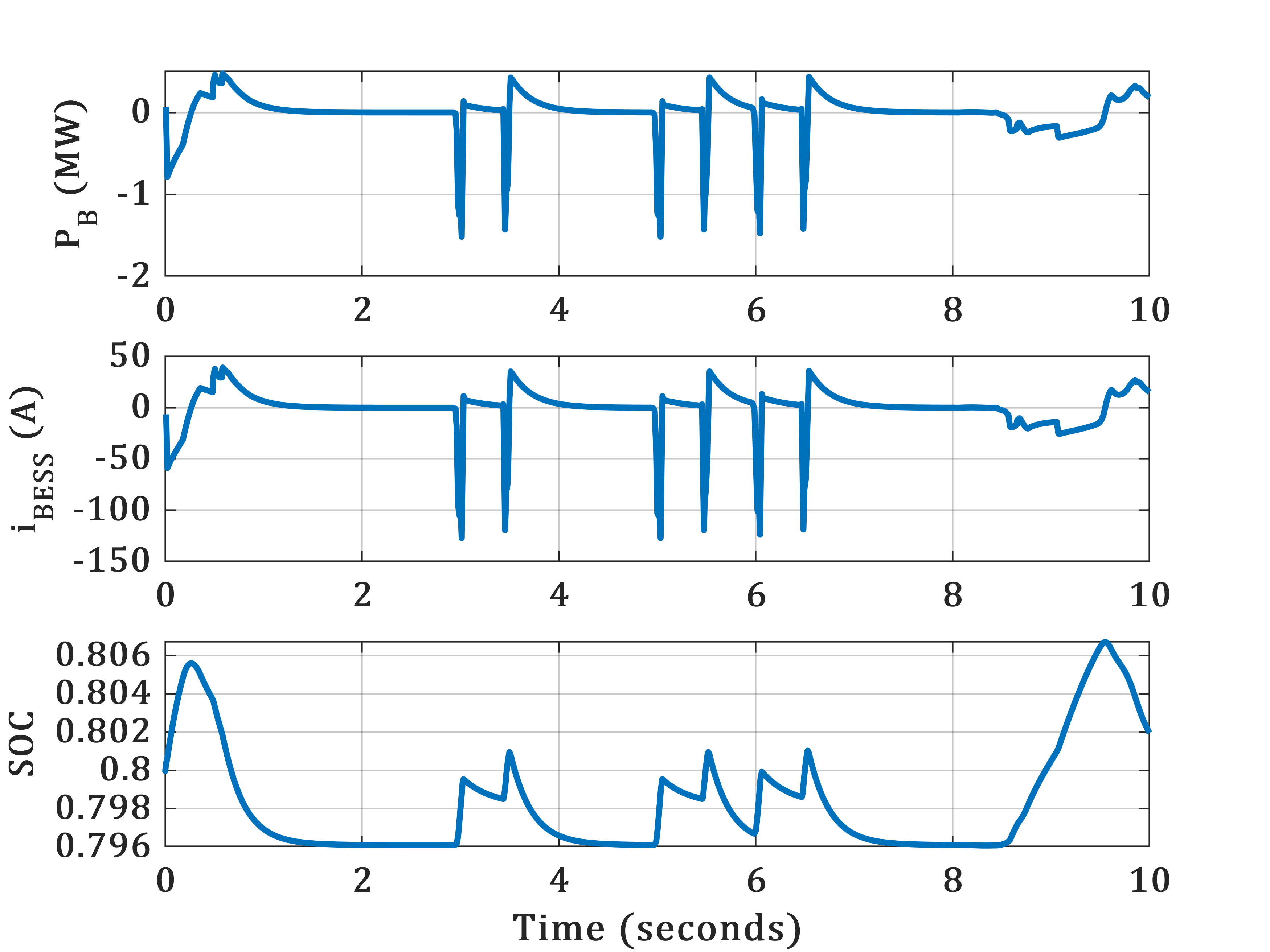} \vspace{-3mm}
 	\caption{(a) Power injected by the ESS (b) ESS real-time current (c) SOC variation for initial SOC of 0.8 and final SOC of 0.8 
     }
 	\label{BESS0808} \vspace{-4mm}
\end{figure}

The system in (\ref{SPS_all}), with generator control in (\ref{SPS_PGM_Cntrl}), battery control in the previous step, parameters in Table \ref{sys_Param}, and the specified load profile is simulated in SIL environment for $10s$. Fig. \ref{Generator08-08} shows the corresponding contributing currents for the generator and the battery ESS. Fig. \ref{BESS0808}c shows the initial and parking SOC corresponding to the processed energy throughout the simulation run. It can be noticed that the MPC enforces the SOC to be as close as possible to target SOC at the end of the horizon. Fig. \ref{Generator08-08}b shows the corresponding feedforward and feedback control of the PGM and Fig. \ref{Vb0808} shows the maintained bus voltage of $12kV$. It can be seen that the PGM is controlled to effectively maintain the main bus voltage levels.   

\subsection{Arbitrary Parking SOC}
In this case the ESS injects power appropriately when required that is in cases where the ramping of generator is not met. It is important to note that the arbitrary value for the final SOC in (\ref{MPC_Const2}) should be carefully choosen based on the ESS size and its workload within the MPC horizon. Otherwise, for example, reaching the final SOC of 0.1 within the optimization horizon will not yield a feasible solution. ESS maximum workload requirements versus high ramp-rate loads are studied in \cite{2019_Bijaieh}.

The system with PGM, DLC, and load profile as the previous subsection is simulated in SIL environment for $10s$. The corresponding power injected by the gen-set and the ESS and their corresponding voltage and current responses and the SOC are presented in Fig. \ref{Generator08-077} and Fig. \ref{BESS08-077}. From Fig. \ref{Vb08-077} it can also be seen that the load bus voltage is regulated at the desired value of 12kV. The combined power injected by the gen-set and ESS versus the power forecast is presented in Fig. \ref{PD08-077}. It can be seen that the power forecast is met by both power injecting components over the target horizon.

\section{CONCLUSIONS}\label{Sec: Conclusion}

This paper presented a framework as a proof of concept to employ a model and load predictive control for effective power and energy management of SPSs for certain applications. It was shown that the generator and the battery systems are controlled according to their specific rated power and ramp-rate limitations and the battery SOC can be managed to park to an arbitrary value. The functionality of the control was verified in a real-time simulation platform. For future work, the demonstrated framework will be employed to specify requirements for events such as undergoing a health-monitoring and degradation identification process and for initial sizing of energy storage elements versus specified operational events while considering the cost and degradation of generators and energy storage devices. Moreover, the results of this work will further be extended for distributed aggregator-based MPC operation in multi-zone SPSs. 



\section{ACKNOWLEDGEMENT}\label{Sec: ACKNOWLEDGEMENT}
This material is based upon research supported by, or in part by, the U.S. Office of Naval Research (ONR) under award number N00014-16-1-2956.









\medskip
\printbibliography

\end{document}